%% file: portrait.tex
%\magnification=1200
\input fig.tex

\def\endofproof{ \smallskip \hskip 4in $\#\#\#$ }
\bigskip
\bigskip
\centerline{\bf On the realization of Fixed Point Portraits.}
\centerline{\bf (An addendum to ``Fixed Point Portraits" by Goldberg and Milnor)}
\bigskip
\centerline{Alfredo Poirier}
\centerline{Mathematics Department}
\centerline{SUNY StonyBrook, NY 11790}
\bigskip
\bigskip
\centerline{\bf Abstract.}
\bigskip
We establish that every formal critical portrait (as defined in [GM]), 
can be realized by a postcritically finite polynomial.
\bigskip
\centerline{\bf 1. Preliminaries.} 
\bigskip
{\bf 1.1.} 
Let $P$ be a polynomial of degree $d \ge 2$ with connected Julia set $J(P)$. 
The {\it rational type $T(z)$} of a fixed point $z$, is by definition the set of all angles  of (rational) external rays which land at $z$. 
The {\it fixed point portrait of $P$} 
is the collection ${\cal T}(P)=\{T_1,\dots,T_k\}$ consisting of all rational types $T_i \ne \emptyset$ of its fixed points. 
\medskip
In the work [GM], Goldberg and Milnor gave combinatorial conditions on the family ${\cal T}(P)$ 
and conjectured that those conditions where also sufficient. 
The purpose of this note is to prove this conjecture. 
\medskip
{\bf 1.2 Rational Rotation Sets} (See also [G].) 
We start by parametrizing the unit circle $S^1={\bf R/Z}$ by the interval $[0,1)$. 
Let $d \ge 2$ and consider the {\it d-}fold covering map\break
$f_d:\theta \mapsto d\theta \hbox{ {\it(mod 1)}}$. 
We will adopt the convention throughout that an indexed subset $\Theta=\{\theta_0,\dots,\theta_{n-1}\}$ of $\bf R/Z$ satisfies 
$0 \le \theta_0 < \dots < \theta_{n-1} <1$. 
\smallskip
{\bf Definition.} 
A finite subset $\Theta=\{\theta_0,\dots,\theta_{n-1}\}$ of $\bf R/Z$ is a {\it degree d-rotation set} if there exists a positive integer $m$ so that 
$f_d(\theta_i)=\theta_{i+m \hbox{ {\it(mod 1)}}}$ for $i=0,\dots,n-1$. 
Note that in this definition $m$ and $n$ need not be relatively prime. 
The ratio $m/n$ {\it (mod 1)} is called the {\it rotation number.}
\medskip
{\bf Theorem.} (See [G, Theorem 7].) 
{\it Let $\Theta,\Theta'$ be degree d rotation sets with the same rotation number $m/n$. 
Then $\Theta=\Theta'$ if and only if for all $i=0,\dots,d-2$ 
$$
\#\big(\Theta \cap [{i \over d-1},{i+1 \over d-1})\big)=\#\big(\Theta' \cap [{i \over d-1},{i+1 \over d-1})\big).
$$
}
In other words a $d-$rotation set is uniquely determined by, 
its rotation number, 
its cardinality, and 
the relative position of its elements with respect to the $d-1$ roots of unity. 
\medskip
{\bf 1.3 Unlinked sets.} 
We will say that two subsets $T$ and $T'$ of the circle $\bf R/Z$ are unlinked if 
they are contained in disjoint connected subsets of $\bf R/Z$, 
or equivalently, 
if $T'$ is contained in just one connected component of the complement ${\bf R/Z}-T$. 
(In particular $T$ and $T'$ must be disjoint.) 
If we identify $\bf R/Z$ with the boundary of the unit disk, 
an equivalent condition would be that the convex closures of these sets are pairwise disjoint. 
As an example, 
if $T$ and $T'$ are the types for any two distinct fixed points of $P$, then evidently $T$ and $T'$ are unlinked. 
\medskip
{\bf 1.4 } 
We fix an integer $d \ge 2$ and a family ${\cal T}=\{T_1,\dots,T_k\}$ to which we impose the following conditions 
\smallskip 
{\it P1. Each $T_j$ is a degree d-rotation set.}

{\it P2. The $T_j$ are disjoint and pairwise unlinked.}

{\it P3. The union of those $T_j$ which have rotation number zero is precisely equal to the set $\{0,{1 \over d-1},\dots,{d-2 \over d-1}\}$ consisting of all angles which are fixed by $f_d$.}

{\it P4. Each pair $T_i \ne T_j$ with non-zero rotation number is separated by at least one $T_\ell$ with zero ration number. 
That is, $T_i$ and $T_j$ must belong to different connected components of the complement ${\bf R/Z}-T_\ell$.} 
\medskip
The importance of the above conditions is shown by the following Theorem proved by Goldberg and Milnor. 
\medskip
{\bf Theorem.} ([GM, Theorem 3.8]) 
{\it If ${\cal T}(P)=\{T_1,\dots,T_k\}$ is the fixed point portrait for some polynomial $P$ with connected Julia set $J(P)$, 
then conditions P1-P4 above are satisfied.}  
\medskip
The main result of this note is the sufficiency of these conditions. 
\bigskip
{\bf Theorem A.} 
{\it Given a family ${\cal T}=\{T_1,\dots,T_k\}$, satisfying conditions P1-P4 above, 
there is a postcritically finite polynomial $P$ such that ${\cal T}(P)={\cal T}$.} 
\medskip
Goldberg and Milnor proved this theorem only for some special cases.  
Our proof is based on the construction of a unique smallest abstract Hubbard tree which realizes the given fixed point portrait. 
(For a different approach, based on Thurston laminations see [HJ].) 
\bigskip
{\bf 1.5 Abstract Hubbard Trees.} 
By an {\it (angled) tree H} will be meant a finite connected acyclic $m$-dimensional simplicial complex ($m=0,1$), 
together with a function 
$\ell,\ell' \mapsto \angle(\ell,\ell')=\angle_v(\ell,\ell') \in {\bf Q/Z}$ 
which assigns a rational modulo 1 to each pair of edges $\ell,\ell'$ which meet at a common vertex $v$. 
This angle $\angle(\ell,\ell')$ should be skew-symmetric, with $\angle(\ell,\ell')=0$ if and only if $\ell=\ell'$, 
and with $\angle_v(\ell,\ell'')=\angle_v(\ell,\ell')+\angle_v(\ell',\ell'')$ whenever three edges are incident at a vertex $v$. 
Such an angle function determines a preferred isotopy class of embeddings of $H$ into ${\bf C}$. 
\medskip
Let $V$ be the set of vertices. 
We specify a mapping $\tau:V \to V$and call it the {\it vertex dynamics}, 
and require that $\tau(v) \ne \tau(v')$ whenever $v$ and $v'$ are endpoints of a common edge $\ell$. 
We consider also a {\it local degree function} $\delta:V \to {\bf Z}$ which assigns an integer $\delta(v) \ge 1$ to each vertex $v \in V$. 
We require that $d(\delta)=1+\sum_{v \in V}(\delta(v)-1)$ be greater that 1. 
By definition a vertex {\it v is critical} if $\delta(v)>1$ and 
{\it non-critical} otherwise. 
The {\it critical set} $\Omega(\delta)=\{v \in V:\hbox{\it v is critical}\}$ is thus not empty. 
\bigskip
The maps $\tau$ and $\delta$ must be related in the following way. 
Extend $\tau$ to a map $\tau:H \to H$ which carries each edge homeomorphically onto the shortest path joining the images of its endpoints. 
We require then that 
$
\angle_{\tau(v)}(\tau(\ell),\tau(\ell'))=\delta(v)\angle_v(\ell,\ell')
$
whenever $\ell,\ell'$
are incident at $v$ (in this case $\tau(\ell)$ and $\tau(\ell')$ are incident at the vertex $\tau(v)$ where the angle is measured). 
\medskip
A vertex $v$ is {\it periodic} if for some $n>0$, $\tau^{\circ n}(v)=v$. 
Given $W \subset V$, we define its orbit ${\cal O}(W)=\cup_{n=0}^\infty \tau^{\circ n}(W)$. 
The orbit of a periodic critical point is a {\it critical cycle}. 
We say that a vertex $v$ is of {\it Fatou type or a Fatou vertex} if it eventually maps into a critical cycle. 
Otherwise, if it eventually maps to a non critical cycle, it is of {\it Julia type or a Julia vertex}.    
\medskip
We define the distance $d_{H}(v,v')$ between vertices in $H$ as the number of edges in a shortest path $\gamma$ between $v$ and $v'$. 
We say that $(H,V,\tau,\delta)$ is {\it expanding} if the following condition is satisfied. 
For any edge $\ell$ whose end points $v,v'$ are Julia vertices, 
there is an $n \ge 1$ such that $d_{H}(\tau^{\circ n}(v),\tau^{\circ n}(v'))>1$. 
Note that angles are not needed in this definition. 
\medskip
The angles at Julia vertices are rather artificial, 
so we normalize them as follows. 
If $m$ edges $\ell_1,\dots,\ell_m$ meet at a periodic Julia vertex $v$, 
then we assume that the angles $\angle_v(\ell_i,\ell_j)$ are all multiples of $1/m$. 
(It follows that the angles at a periodic Julia vertex convey no information beyond the cyclic order of these $m$ incident edges.) 
\medskip
{\bf Definition.} 
By an {\it abstract Hubbard tree} we mean an angled tree 
${\bf H}=((H,V,\tau,\delta),\angle)$ such that the angles at any periodic Julia vertex where $m$ edges meet are multiples of $1/m$. 
We define isomorphism between abstract Hubbard Trees in an obvious way.
\bigskip
Douady and Hubbard showed in [DH] that 
a postcritically finite polynomial $P$ and a finite invariant set $M$ containing the critical set $\Omega(P)$ of $P$ 
naturally defines an abstract Hubbard tree ${\bf H}_{P,M}$. 
To define the angle function we note the following facts. 
Near a Fatou vertex the edges of the tree are by definition segments of constant argument in the B\"ottcher coordinate, 
we define the angle between two such edges as the difference in such coordinates. 
For a (periodic or preperiodic) Julia set point $v$, $J(P)-\{v\}$ consists of a finite number (say m) of components. 
We define the `angle' between consecutive components around $v$ to be $1/m$. 
As each edge in the tree correspond locally to one of these components, 
we have an angle function between them. 
(This procedure is well defined and compatible with the definition above, see [P].) 
It is easy to prove that this abstract Hubbard tree is expanding (see [P]). 
\medskip
The main result for Hubbard trees is the following. 
\medskip
{\bf Theorem 1.} (See [P].)  
{\it Let ${\bf H}$ be an abstract Hubbard tree. 
Then there is a postcritically finite polynomial P and an invariant set $M \supset \Omega(P)$ such that ${\bf H}_{P,M} \cong {\bf H}$ 
if and only if $\bf H$ is expanding. 
Furthermore, P is unique up to affine conjugation.} 
\medskip
This abstract Hubbard tree also gives information about external rays as the following theorem essentially due to Douady and Hubbard shows 
(see [DH, Chap VII] or [P]).
\smallskip
{\bf Theorem 2.}
{\it The number of rays which land at a periodic Julia vertex is equal to 
the number of incident edges of the tree, 
and in fact, 
there is exactly one ray landing between each pair of consecutive edges. 
Furthermore, the ray which lands at $v$ between $\ell$ and $\ell'$ maps to the ray 
which lands at $\tau(v)$ between $\tau(\ell)$ and $\tau(\ell')$.}
\bigskip
\bigskip
\centerline{\bf 2. Proof of Theorem A.}
\bigskip
We identify $\bf R/Q$ with $\partial D$ via the exponential map $e(\theta)=e^{2\pi i \theta}$. 
For each element $T_j=\{\theta_1,\dots,\theta_n\}$ consider the baricenter $v(T_j)$ of all elements of $T_j$. 
In other words define 
$$
v(T_j)={1 \over n}\sum_{i=1}^n e(\theta_i). 
$$
Next we join each element $e(\theta) \in T_j$ to $v(T_j)$ by a straight segment 
(these segments will not be part of our tree). 
This construction clearly divides the closed unit disk into a finite number of components or {\it regions} (see condition P2 in $\S1$ and compare Figure 1). 
\medskip 
For each of these regions $U_i$ define {\it the critical capacity} $CC(U_i)$
as the number of vertices $v(T_j)$ with zero rotation number belonging to the boundary of $U_i$. 
Clearly the sum of critical capacities must be equal to $d-1$ 
(this is an easy induction using condition $P3$).   
\medskip 
Insert inside each region $U_i$ a vertex $w(U_i)$, which we join to every vertex $v(T_j)$ in the boundary of $U_i$. 
The union of these joining edges together with the vertices $w(U_i)$ and $v(T_j)$ will form the required topological tree (compare with Figure 2). 
\medskip
We proceed now to construct the local degree and angle functions of the tree. 
For every vertex $w(U_i)$ we define its degree $\delta(w(U_i))=CC(U_i)+1$.  
For the vertices $v(T_j)$ define $\delta(v(T_j))=1$. 
At every vertex where $m \ge 1$ edges come together we define the angle between consecutive edges to be $1/m$. 
At a vertex $w(U_i)$ this number equals to the number of vertices $v(T_j) \in \partial U_i$, 
while at a vertex $v(T_j)$ it is equal to the number of elements in $T_j$. 
\medskip
To define the vertex dynamics we consider first 
those vertices $w(U_i)$ for which the region $U_i$ has a (necessarily unique by condition $P4$) vertex $v(T_l)$ with non zero rotation number on its boundary. 
Then, there are exactly two elements $\theta,\theta' \in T_l$ such that 
$e(\theta),e(\theta')\in \partial U_i$. 
In fact, we can order these two elements so that 
$e(\theta+\epsilon),e(\theta'-\epsilon)\in \partial U_i$ for small $\epsilon>0$. 
Then there is a unique $U_j$ $(\ne U_i)$ such that 
$e(d\theta+\epsilon),e(d\theta'-\epsilon)\in \partial U_j$, 
and we define $\tau(w(U_i))=w(U_j)$. 
(Note that in this way $w(U_i)$ has the same period as any $\theta \in T_l$ under multiplication by $d$ modulo $1$.) 
For all other vertices $v$ define $\tau(v)=v$. 
\medskip
Note that by construction all $v(T_l)$ are of Julia type (fixed and non critical), 
while all $w(U_l)$ are of Fatou type. 
Also, 
between any two different $v(T_i)$ and $v(T_j)$ there is a vertex $w(U_l)$, 
and so, the expanding condition is trivially satisfied. 
Furthermore, at non fixed Fatou vertices 
the angle between consecutive edges is $1/\delta$, 
while at fixed Fatou vertices it is $1/(\delta-1)$. 
From this it is easy to see that the angle condition is satisfied. 
\smallskip
Thus, there is a unique (up to affine conjugation) polynomial of degree d which realizes this abstract Hubbard Tree. 
We must still verify that this polynomial (or tree) has the required fixed point portrait. 
We begin by locating all fixed points. 
\medskip
{\bf Lemma.} 
{\it The abstract Hubbard tree constructed above has exactly d fixed points.} 
\medskip
{\bf Proof.} 
This is just a matter of counting. 
Let $k$ be the number of rotation sets,\break
and let $\ell=\sum \#\{T_i \hbox{ : } T_i$ has non zero rotation number$\}$. 
Note that by condition $P3$\break  
$d-1=\sum \#\{T_i \hbox{ : } T_i$ has zero rotation number$\}$.  
By induction it is easy to see that there are exactly $\ell+d-k=(1+\sum(\#(T_j)-1))$ regions $U_l$. 
There are $k$ Julia fixed points (as many as rotation sets). 
By construction there are $\ell$ regions without an interior fixed point,  
so there are $d-k$ Fatou fixed points. 
\endofproof
\medskip
To verify that this tree has the required fixed point portrait we use theorem 2 in section 1.5. 
We first note that by construction all edges incident at a vertex $v(T_l)$ where $T_l$ has rotation number zero are `fixed'. 
Thus by theorem 2, 
only fixed rays land there, 
and every fixed ray must land at one of those points. 
Now suppose $0 \in T_\ell$, then 
by construction (if $T_\ell$ is not a singleton) there is a segment (the one which is not part of the tree!) joining $e(0)$ to $v(T_\ell)$. 
As topologically this segment is located between two consecutive fixed edges of the tree, 
it corresponds to a fixed ray of the polynomial map (this is true even if $T_\ell$ is a singleton). 
After affine conjugation if necessary we assume that this is the zero ray. 
Now, if we walk counterclockwise around the tree, 
it follows from conditions $P2$ and $P3$ that given $T_j$ with zero rotation number, 
all rays with argument in $T_j$ land at $v(T_j)$.  
\smallskip
If $T_l$ has rotation number $m/n$, clearly 
the set of rays which land at $v(T_l)$ has also rotation number $m/n$. 
The result then follows from Theorem 1.2, 
which asserts that a $m/n$ rotation set is uniquely determined by the relative position of its elements respect to ${i \over d-1}$, $i=0,\dots,d-2$. 
\bigskip
\centerline{\bf 3. An example.} 
\bigskip
\bigskip
We will illustrate the proof of the above theorem by taking the degree 5 fixed point portrait determined by 
$T_1=\{0,3/4\}$, $T_2=\{1/8,5/8\}$, $T_3=\{1/4\}$, $T_4=\{1/2\}$, 
with rotation numbers $0$, $1/2$, $0$, $0$ respectively. 
The Julia set of the actual polynomial is shown in figure 3. 
\bigskip
\centerline{\it \hskip 0.35in Step 1: \hskip 2.6in Step 2:}
\centerline{
\hbox{
\insertRaster 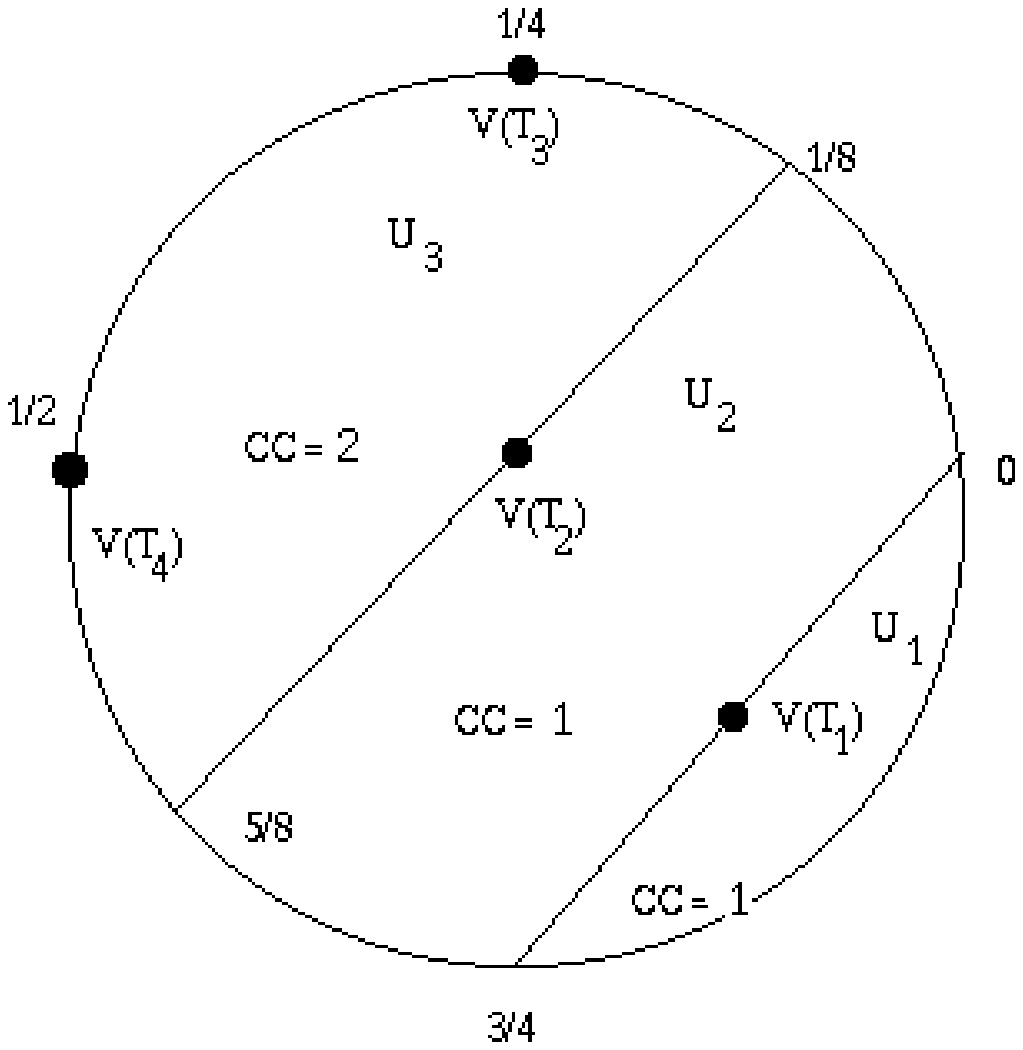 pixels 344 by 319 scaled 500
}
\hskip 0.75in
\hbox{
\insertRaster 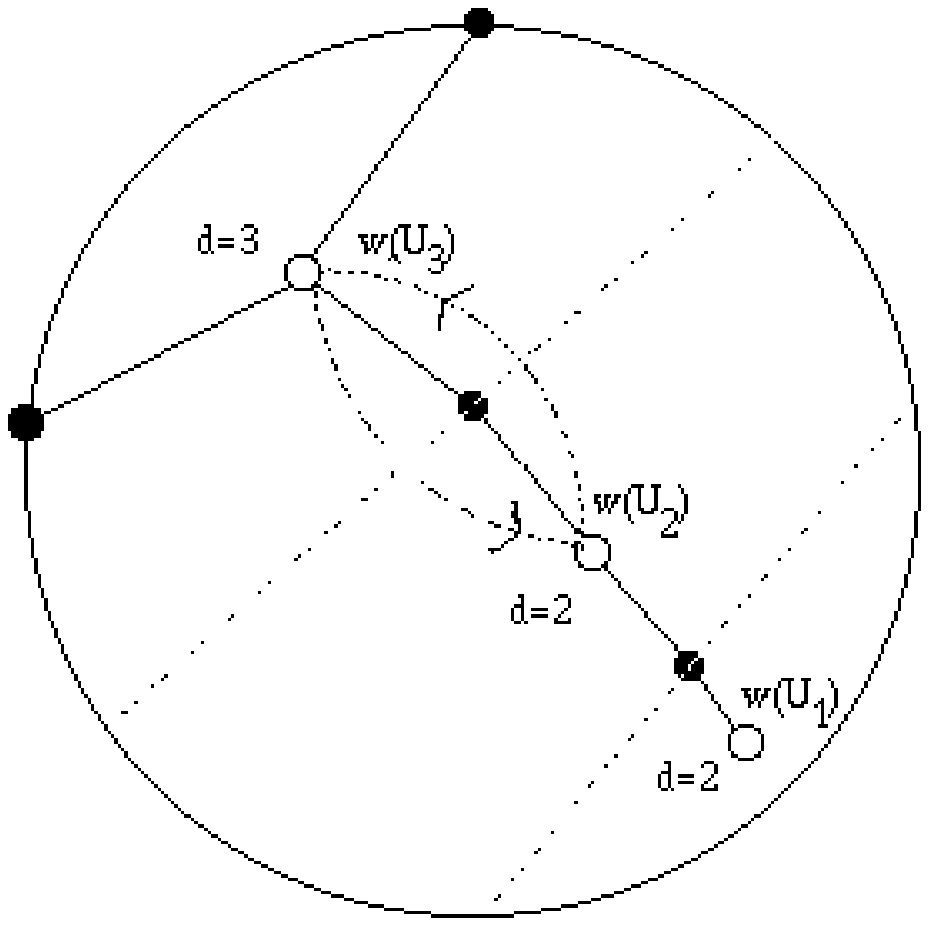 pixels 279 by 277 scaled 600
}}
\medskip
\centerline{\it \hskip 0.25in Figure 1 \hskip 2.5in Figure 2}
\centerline{\it \hskip 3in all vertices are fixed except}
\centerline{\hskip 3in {\it $w(U_2)$, $w(U_3)$ which are interchanged.}}
\vfill\eject
\centerline{
\insertRaster 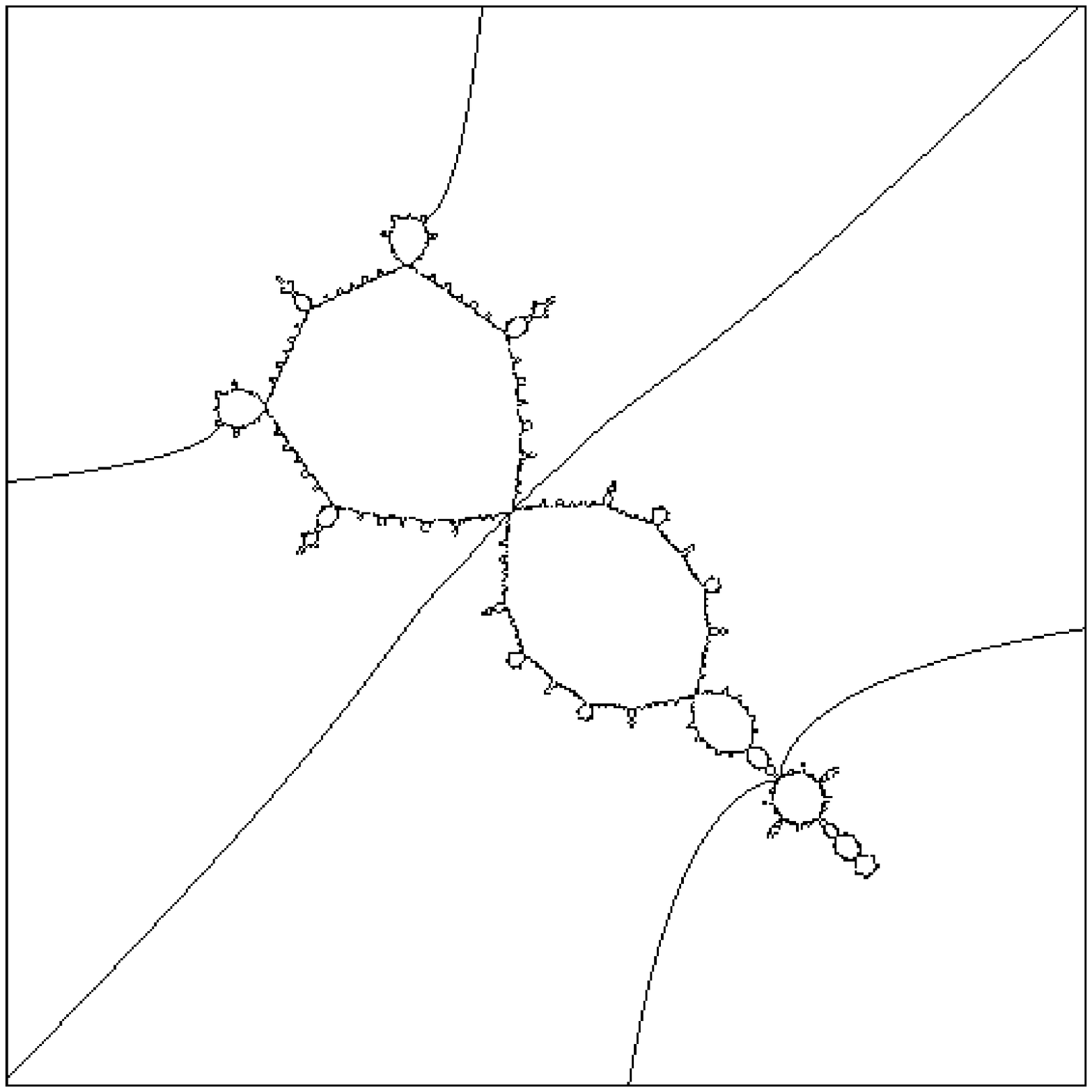 pixels 480 by 480 scaled 400
}
\medskip
\centerline{\it  \hskip 0.175in Figure 3}
\bigskip
\centerline{\it Julia set of the polynomial $P(z)=z^5+Az^3+Bz^2+Cz+D$, where}
\centerline{\it $A \approx 2.714670827i$, $B \approx 0.693957313(1+i)$, $C \approx -1.608651885$}
\centerline{\it $D \approx -0.355745059(1-i)$. }
\centerline{\it The rays $0,{1 \over 8},{1 \over 4},{1 \over 2},{5 \over 8},{3 \over 4}$ shown.}
\bigskip
{\bf Acknowledgement.} 
We will like to thank John Milnor for helpful conversations. 
The inclusion of the example as well as several other remarks were suggested by him. 
Figure 3 was constructed using a program of Milnor. 
Also, we want to thank the Geometry Center, University of Minnesota
and Universidad Cat\'olica del Per\'u 
for their material support. 
\bigskip
\bigskip
\centerline{\bf References.}
\bigskip
\bigskip
[DH] A.Douady and J.Hubbard, \'Etude dynamique des polyn\^omes complexes, part I; 
Publ Math. Orsay 1984-1985. 
\smallskip
[G] L.Goldberg, Rotation Subsets; Preprint \#1990/14 IMS SUNY@StonyBrook.
\smallskip
[GM] L.Goldberg and J.Milnor, Fixed Point Portraits; Preprint \#1990/14
IMS\break SUNY@StonyBrook.
\smallskip
[HJ] S.Hu and Y.Jiang, Toward Topological Classification of Critically Finite Polynomials; (in preparation).
\smallskip
[P] A.Poirier, On Postcritically Finite Polynomials; Thesis, SUNY@StonyBrook, 1992 (to appear). 
\bigskip
\bigskip
Typeset in $\TeX$ (November 25, 1991) 
\end

%% file: fig.tex
% This file is supposed to imitate the -ras.tex- file which contained
% the macros -insertRaster- and -RasterBox- for including SUN raster files
% in TeX.  But here we assume that the original raster files have been
% converted to PostScript, so we use -psfig- rather than -sunbitmap-.
% There was a bug in the way -ras.tex- scaled these bitmaps, and this file
% attempts to preserve this bug so that the resulting figures are identical
% to the raster figures in the original paper.

\input psfig
%\magnification=1200

\newif\ifboxfigure      % set to true if you want the figure boxed
\boxfigurefalse

\def\BoxIt#1#2{%        % put a box around #1, leaving a gap of #2
	\vbox{\hrule
	\hbox{\vrule\kern#2\vbox{\kern#2#1\kern#2}\kern#2\vrule}
		   \hrule}}

\def\insertRaster #1 pixels #2  by #3 scaled #4 {
%  #1 is the raster file
%  #2 is the width of the picture in pixels
%  #3 is the height of the picture in pixels
%  #4 is scaling 500 means half size 2000 means double size 1000 is actual  
			        %  \medskip
			 \hbox {%     to \hsize{%

			 \hss
			 \RasterBox {#1} {#2} {#3} {#4}
			 \hss
			 }%
}

\def\RasterBox #1 #2 #3 #4{

% dimen0 and dimen1 are the dimensions of the hbox which are setup in
% the same way as -ras.tex- unfortunately psfig includes images in a 
% different manner from sunbitmap, so it must be scaled using different
% values, hence dimen2 and dimen3 and MAGIC number 65/72

\dimen5=65pt
\divide\dimen5 by 72

\dimen0=#2\dimen5
\divide\dimen0 by 1000
\dimen1=#3\dimen5
\divide\dimen1 by 1000
\dimen2=#3\dimen5
\divide\dimen2 by 1000
\dimen3=#2\dimen5
\divide\dimen3 by 1000

\setbox4=\hbox to #4\dimen0{
 \vbox to #4\dimen1{
 \vss
 \psfig{figure=#1,height=#4\dimen2,width=#4\dimen3}
 }
 \hss
 }
 \ifboxfigure\BoxIt{\box4}{0pt}
 \else\box4
 \fi
 }

% \insertRaster towrrab.pic pixels 500 by 500 scaled 400